\documentclass[10pt]{article} 
\usepackage{amsmath, amssymb, amsthm, mathrsfs}
\usepackage{enumerate}
\usepackage{colordvi}
\usepackage{tikz}
\usepackage{graphicx}
\usepackage{hyperref}
\newtheorem{claim}{Claim}[section]
\newtheorem{theorem}[claim]{Theorem}
\newtheorem{lemma}[claim]{Lemma}

\newtheorem{corollary}[claim]{Corollary}
\usepackage{color}

\definecolor{Myred}{cmyk}{0.0,1.0,1.0,0.00}

\newtheorem{definition}{Definition}

\begin{document}
%\begin{center}
%{\Large{\textbf{Spectral convergence of Neumann Laplacian in the perforated domain}}}
\title{Neumann Laplacian in a perturbed domain}
%\bigskip
\author
{
Diana Barseghyan$^{1,2}$, Baruch Schneider$^{1}$ and Ly Hong Hai$^{1}$
}
\date{\small $^{1}$ Department of Mathematics, University of Ostrava,  30.dubna 22, Ostrava 70103, Czech Republic\\ 
$^{2}$ Nuclear Physics Institute, Czech Academy of Sciences, 250 68 Rez, \\ Czech Republic
\\
E-mails:\, diana.barseghyan@osu.cz, baruch.schneider@osu.cz, P18180@student.osu.cz
}
\maketitle

\begin{abstract}
We consider a domain with a small compact set of zero Lebesgue measure of  removed. Our main result concerns the spectrum of the Neumann Laplacian defined on such domain. We prove that the
spectrum of the Laplacian converges in the Hausdorff distance sense to the spectrum of the Laplacian defined on the unperturbed domain.
\end{abstract}
\bigskip


%%%%%%%%%%%%%%%%%%%%%%

\section{Introduction}
\label{s:intro}

What  kind  of  spectral convergence can we expect for the Laplace operator under  perturbations of the domain such as removing small holes?  
It is a common expectation that small perturbations of the physical situation lead only to a small change of the spectrum. In the case of domain perturbations this is largely true for Dirichlet boundary conditions while the Neumann case is more delicate.

Such questions received already quite a lot of answers starting from the seminal work of Rauch and Taylor
concerning the spectrum of the Laplace operator of the domains with holes \cite{RT75}. 

In Neumann case even small perturbations may cause abrupt change of the spectrum. For example, such an effect is observed when the hole has a "split-ring" geometry \cite{S15}. The split ring (even being very small) may produce additional eigenvalues having nothing in common with the eigenvalues of the Neumann Laplacian on the unperturbed domain.

In \cite{O83} Ozawa have studied the asymptotic behavior of the eigenvalues of Laplace operator on the domain with small spherical obstacles imposing the Neumann condition at their boundary and the Dirichlet condition at the rest part of the boundary.

Maz’ya, Nazarov and Plamenewskii, see \cite[Ch.9, vol.I]{MNP00}, have  considered the Laplace operator on the domain with obstacles, imposing the Dirichlet condition on their boundary and have proved the validity of a complete asymptotic expansion for the eigenvalues.

In \cite{BN98}, \cite{CP20}, \cite{D03}, \cite{KP17}  authors have considered the Dirichlet Laplacians on Euclidean domains or manifolds with holes and studied the problems of the resolvent convergence.

The problems with Neumann obstacles having more general geometry appeared in \cite{CP20}.
The authors required the hole to satisfy the so called "uniform extension property". This is requirement, which means that $\mathcal{H}^1$- functions on the domain with a hole can be extended to $\mathcal{H}^1$ function on the unperturbed domain and the norm of this extension operator does not depend on the hole diameter. In this case the authors established the spectral convergence.

Another related paper, which makes use imposing Neumann boundary conditions on the boundary of hole which has satisfies the suitable geometrical assumptions to consider  the upper and lower estimates for the ground state eigenvalue, have been extensively studied by Hempel, see for instance \cite{H06}.

There is one case, which was not considered in \cite{CP20} (and other similar papers). These are holes with zero Lebesgue measure (e.g. an interval or a piece of a curve). Evidently, for the holes with Lebesgue measure zero such an extension required in \cite{CP20} is not possible. 

We will be interested in a two-dimensional bounded domain with a single hole $K_\varepsilon$ (for a fixed  parameter $\varepsilon$) having zero Lebesgue measure. The main purpose of this paper is to prove the spectral convergence of the Neumann Laplacian  on  $\Omega_{K_\varepsilon}$ as $\varepsilon\to 0$ in terms of the Hausdorff distance under some additional assumptions on the geometry of $K_\varepsilon$.

The outline of the paper is as follows. It consists of four sections besides the introduction.
In the second section we present the main  tools  of  the spectral convergence  of  operators  on varying Hilbert spaces. In the third section we present the main results, and consider a general theorem, namely Theorem \ref{Main}. We will use Theorem  \ref{Main} in the proof of Theorem \ref{Cons.}, and where consider the spectral convergence for the Laplacians on $\Omega$ and $\Omega_{K_\varepsilon}$. In the fourth and last section we prove our results to which we already alluded.

%%%%%%%%%%%%%%%%%%%%%%%%%%%%%%%%%%%%%%%%
\section{Main  tool  of the spectral convergence  of  operators on varying Hilbert spaces}
\setcounter{equation}{0}

We begin this section by reviewing some basic facts which assures a spectral convergence for two operators having the different domains. For more information we refer the reader to  \cite{P06}. 

To a Hilbert space $H$ with inner product $(\cdot, \cdot)$ and norm $\|\cdot\|$ together with a non-negative, unbounded, operator $A$, we associate the scale of Hilbert spaces
\[
H_k:=\mathrm{Dom}((A+I)^{k/2}),\quad \|u\|_k:=\|(A+ I)^{k/2}u\|,\,\,k\ge 0,
\]
where $I$ is the identity operator.

We think of $(H', A')$ being some perturbation of $(H, A)$ and want to lessen the assumption such that the spectral properties are not the same but still are close.

\begin{definition}(see \cite{P06})\label{Post} 
Suppose we have linear operators
\begin{gather*} J: H \longrightarrow H', \quad\quad J_1: H_1\longrightarrow H_1'\\
J': H'\longrightarrow H,\quad\quad J_1': H_1' \longrightarrow H_1. \end{gather*}
Let $\delta> 0$ and $k\ge 1$.
We say that $(H, A)$ and $(H', A')$ are $\delta$-close of order $k$ iff the following conditions are fulfilled:

\begin{gather}\label{1} \|J f- J_1 f\|_0\le \delta \|f\|_1, \\\label{2}  |(J f, u)-(f, J' u)|\le \delta \|f\|_0 \|u\|_0,\\\label{3}  \|u-J J' u\|_0\le \delta \|u\|_1,\\\label{4} \quad \|Jf\|_0\leq 2\|f\|_0,\,\,\, \|J'u\|_0\le 2\|u\|_0,
\\\label{5'} \|(f- J'J f)\|_0\le \delta \|f\|_1,\\
\label{6} \| J' u- J_1'u\|_0\le \delta \|u\|_1\\\label{7}  |a(f, J_1' u)- a'(J_1 f, u)|\le \delta \|f\|_k \|u\|_1,\end{gather}
for all $f, u$ in the appropriate spaces. Here, $a$ and $a'$ denote the sesquilinear forms associated to $A$ and $A'$.
\end{definition}
We denote by $d_{\mathrm{Haussdorff}}(A, B)$ the Hausdorff distance for subsets $A, B\subset\mathbb{R}$
\begin{equation}\label{Haus.}
d_{\mathrm{Haussdorff}}(A, B):=\mathrm{max} \left\{\underset{a\in A} {\mathrm{sup}}\,d(a, B),\,\underset{b \in B} {\mathrm{sup}}\, d(b, A)\right\},
\end{equation}
where $d(a, B):= \mathrm{inf}_{b\in B} |a - b|$. We set
\begin{equation}\label{Hausdorff2}\overline{d} (A, B):= d_{\mathrm{Hausdorff}}\left((A+I)^{-1}, (B+I)^{-1}\right)\end{equation}
for closed subsets of $[0, \infty)$. For the next result, which originates with the  work of O. Post \cite{P06} we have the following spectral convergence theorem in terms of the distance $\overline{d}$.
\begin{theorem}\cite{P06}\label{POSTtm}
There exists $\eta(\delta)>0$ with $\eta(\delta)\to 0$ as $\delta\to 0$ such that
\begin{equation}\label{Post}\overline{d}(\sigma_\bullet(A), \, \sigma_\bullet (A'))\le \eta(\delta)\end{equation}
for all pairs of non-negative operators and Hilbert spaces $(H, A)$ and $(H', A')$ which are $\delta$-close. Here, $\sigma_\bullet (A)$ denotes either the entire spectrum, the essential or the discrete spectrum of $A$. Furthermore, the multiplicity of the discrete spectrum, 
$\sigma_{\mathrm{disc}}$, is preserved, i.e. if $\lambda\in \sigma_{\mathrm{disc}}(A)$ has multiplicity $\mu>0$, then there exist $\mu$ eigenvalues (not necessarily all distinct) of operator $A'$ belonging to interval $(\lambda-\eta(\delta), \lambda+\eta(\delta))$.
\end{theorem}

\section{Main results}
\setcounter{equation}{0}

The starting point is to consider a bounded domain $\Omega\subset \mathbb{R}^2$ and a compact set $K\subset \Omega$ with zero Lebesgue measure (e.g. an interval or a piece of a curve). We denote $\Omega_K:= \Omega\setminus K$. 
The Neumann Laplacian $-\Delta_N^{\Omega_K}$ is defined on the Sobolev space $\mathcal{H}^1(\Omega_K)$ via the quadratic form
\[
\int_\Omega |\nabla u|^2\,d x\,d y,\quad u\in \mathcal{H}^1(\Omega_K).
\]
In case if $K$ is  empty set then unperturbed Neumann Laplacian denoted by $-\Delta_N^\Omega$ which is defined via the same form $\int_\Omega |\nabla u|^2\,d x\,d y$ on $\mathcal{H}^1(\Omega)$.

\bigskip
During our paper we suppose the additional \textbf{property$^*$}: for any $(x_0, y_0)\in \Omega_K$
at least one of the following conditions takes place
\begin{itemize}
\item The line $l(x_0)=\{x= x_0\}\subset \Omega_K$ has no intersection with $K.$   

\item The line $h(y_0)=\{y= y_0\}\subset \Omega_K$ has no intersection with $K.$  
\end{itemize}
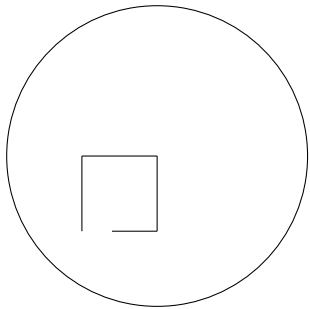
\begin{tikzpicture}
\draw (2,2) -- (2,1);
\draw (1,2) -- (2,2);
\draw (1,2) -- (1,1);
\draw (1.4,1) -- (2,1);
\draw (2,2) circle (2cm);
\end{tikzpicture}

%%%%%%%%%%%%\textbf{PICTURE}
\bigskip

The main result of this section is the following theorem.

\begin{theorem}\label{Main}
Let $\Omega$ be an open bounded domain in $\mathbb{R}^2$ and let $p\in \Omega$ be some fixed point. Suppose that ${\mathbb B}_\varepsilon\subset \Omega$ is a ball with center at $p$ and radius $\varepsilon>0$. Let $K= K(\varepsilon)\subset \mathbb B_\varepsilon$ be a compact set with zero Lebesgue measure (e.g. an interval or a piece of a curve). Moreover, suppose that $\Omega_K$ satisfies property$^*$. Let $-\Delta_\Omega^N$ and $-\Delta_{\Omega_{K_\varepsilon}}^N$ be the Neumann Laplacians defined on $\Omega$
and $\Omega_{K_\varepsilon}$, respectively. Then for small enough $\varepsilon$ the Neumann Laplacians $-\Delta_\Omega^N$ and $-\Delta_{\Omega_{K_\varepsilon}}^N$ are  $ \mathcal{O}\left(\varepsilon^{1/6}\right)$  close of order $2$.
\end{theorem}
As a consequence we formulate the following theorem:
\begin{theorem}\label{Cons.}
Using the above assumptions, let $-\Delta_\Omega^N$ and $-\Delta_{\Omega_{K_\varepsilon}}^N$ be the Neumann Laplacians defined on $\Omega$
and $\Omega_{K_\varepsilon}$, respectively.
Then there exists $\eta(\varepsilon)>0$ with $\eta(\varepsilon)\to 0$ as $\varepsilon\to 0$ such that the following spectral convergence takes place
$$\overline{d} \left(\sigma_{\bullet} \left(-\Delta_{\Omega_{K_\varepsilon}}^N\right), \, \sigma_{\bullet} \left(-\Delta_\Omega^N\right)\right) \le \eta(\varepsilon),$$
where $\overline{d}$ is defined in (\ref{Hausdorff2}) and $\sigma_\bullet (\cdot)$ denotes either the entire spectrum, the essential or the discrete spectrum. Moreover, the multiplicity of the discrete spectrum is preserved.
\end{theorem}
{\bf Proof of Theorem \ref{Cons.}}. The proof follows directly from the inequality (\ref{Post}) and Theorem \ref{Main}.
\qed
\\

\noindent As consequences of Theorem \ref{Cons.}:

\begin{corollary}
Suppose that $-\Delta_\Omega^N$ has purely discrete spectrum denoted by $\lambda_k(\Omega)$ (repeated according to multiplicity). Then the infimum of the essential spectrum of $-\Delta_{\Omega_{K_\varepsilon}}^N$ tends to infinity and there exists $\eta_k(\varepsilon)>0$ with $\eta_k(\varepsilon)\to 0$ as $\varepsilon\to 0$ such that
$$|\lambda_k(\Omega)- \lambda_k(\Omega_{K_\varepsilon})|\le \eta_k(\varepsilon)$$
for small enough $\varepsilon$. Here, $\lambda_k(\Omega_{K_\varepsilon})$ denotes the discrete spectrum of $-\Delta^N_{\Omega_{K_\varepsilon}}$ (below the essential spectrum) repeated according to multiplicity.
\end{corollary}

\begin{corollary} The Hausdorff distance between the spectra of $-\Delta_{\Omega_{K_\varepsilon}}^N$ and $-\Delta_\Omega^N$ converges to zero on any compact interval $[0, \Lambda]$. 
\end{corollary}

We now turn to the proof of Theorem \ref{Main}.

\section{Proof of Theorem\,\ref{Main}}

We proceed in a number of steps.
\bigskip
\\
STEP\,1:\,{\it Construction of the mappings $J, J', J_1, J_1'$}.

\bigskip

It is easy to notice that
$H=H'= L^2(\Omega),\,A= A'=-\Delta$;
$H_1,\,H_1'$ correspond to Sobolev spaces $\mathcal{H}^1(\Omega)$ and $\mathcal{H}^1(\Omega_{K_\varepsilon})$ and
$H_2= \mathrm{Dom}(-\Delta_N^\Omega)$. The norm $\|\cdot\|_0$ corresponds with the
$L^2$ norm and
$$\|u\|_1=(\|u\|^2_0+ \|\nabla u\|_0^2)^{1/2},\quad
\|f\|_2=\|-\Delta f +f\|_0.$$ 
%%%%%%%%%

Since $H= H'$ and $H_1 \subset H_1'$ we choose $J=J'=I$, where $I$ is the identity operator and 
$J_1$ is the restriction operator:  $J_1u= u\restriction_{\Omega_{K_\varepsilon}}$ for $u\in H_1$.

Let us now construct the mapping $J_1': H_1' \to H_1$.
Without loss of generality, assume that the ball $\mathbb B_\varepsilon$ mentioned in Theorem \ref{Main} is centered at the origin. Let $\epsilon\in (\varepsilon, 2\varepsilon)$ be a number to be chosen later and let $\mathbb B_\epsilon \supset \mathbb B_\varepsilon$ be the ball with center again at the origin and radius 
$\epsilon$, $\Omega_\epsilon:=\Omega\backslash \mathbb B_\epsilon$.

We are going to construct mapping $J_1'$ first for smooth functions. For any $v\in  C^\infty(\Omega_{K_\varepsilon})$ we define

\begin{equation*}J_1' v:=\begin{cases} v,\quad \text{on} \quad \Omega_\epsilon, \\ \frac{r}{\epsilon} \tilde{v}(\epsilon, \varphi), \quad \text{on} \quad \mathbb B_\epsilon,\end{cases}\end{equation*}
where $\tilde{v}(r, \varphi)=v(r \cos \varphi, r \sin \varphi)$.

Now let us construct the mapping $J_1' u$ for any $u\in H_1'$.  Employing the approximation method described in \cite[Thm.2, 5.3.2]{E10}, for the fixed sequence $\{\eta_k\}_{k=1}^\infty$ converging to zero we construct the sequence 
$v_{\eta_k}\in C^\infty(\Omega_{K_\varepsilon})$ which satisfies
\begin{equation}\label{converging}
\|u-v_{\eta_k}\|_1=\|u-v_{\eta_k}\|_{\mathcal{H}^1(\Omega_{K_\varepsilon})}< \eta_k\,\|u\|_1.\end{equation}

Let us mention that in view of the inequalities (\ref{B}) and (\ref{J1'}) which will be proved later it follows  for any smooth function $v$
$$\|J_1'v\|_1^2= \int_\Omega |\nabla J_1'v|^2\,d x\,d y+\int_\Omega |J_1'v|^2\,d x\,d y\le \overline{C}(\varepsilon) \|v\|_1^2,$$
where $\overline{C}(\varepsilon)$ is some constant.

Therefore using the completeness of space $H_1$ we are able to define
\begin{equation}\label{approximation}J_1' u= \underset{k\to \infty}{\mathrm{lim}} \,J_1'v_{\eta_k}.
\end{equation}
\\
STEP\,2:\,{\it The conditions (\ref{1})-(\ref{7}) hold for the mappings $J, J', J_1, J_1'$}.
\bigskip
\\
Indeed, we have that the estimates (\ref{1})-(\ref{5'}) are satisfied with $\delta=0$. 
\\\\

We now prove (\ref{6}), i.e. {\it under the assumptions stated in Theorem \ref{Main} inequality (\ref{6}) is satisfied with $\delta=\mathcal{O}(\sqrt{\varepsilon})$ for small enough $\varepsilon$.}
\\

In view of our construction we have 
\begin{equation}\label{start 6}
\|J' u- J_1'u\|_0^2=\int_\Omega |u- J_1' u|^2\,d x\,d y.
\end{equation}
Let us estimate the right-hand side of (\ref{start 6}). By virtue of (\ref{converging}) and (\ref{approximation}) it is enough to prove (\ref{start 6}) for $u\in C^\infty(\Omega_{K_\varepsilon})$.

Taking into account the relation $J_1' u= u$ on $\Omega_\epsilon$ and using (\ref{approximation}) the above bound implies
\begin{gather}
\nonumber
\int_\Omega |u- J_1' u|^2\,d x\,d y= \int_{\mathbb B_\varepsilon}|u- J_1'u|^2\, d x\,d y\\\label{ineq.}\le 2\int_{\mathbb B_\epsilon}|u|^2\,d x\,d y+ 2\int_{\mathbb B_\epsilon}|J_1'u|^2\, d x\,d y.
\end{gather}

To estimate the first integral of the right-hand side of (\ref{ineq.}) we need the following lemma (the proof is given in the Appendix):

\begin{lemma}\label{curve}
Let $\omega\subset \mathbb{R}^2$ be an open domain satisfying \textbf{property$^*$}.
Then for any $u\in \mathcal{H}^1(\omega)$ and for almost any
 $(x_0, y_0)\in \omega$ the following statement takes place: 
\begin{eqnarray}\label{curve1}
|u(x_0, y_0)|\le C \left(\int_{l(x_0)}\left|\frac{\partial u}{\partial y}(x_0, z)\right|^2\,d z+ \int_{l(x_0)}|u(x_0, z)|^2\,d z\right)^{1/2},\\\nonumber  \text{if} \quad (x_0, y_0) \quad\text{satisfies first condition of \textbf{property$^*$}} 
\\\label{curve2}
|u(x_0, y_0)|\le C \left(\int_{h(y_0)}\left|\frac{\partial u}{\partial z}(z, y_0)\right|^2+ \int_{h(y_0)}|u(z, y_0)|^2\,d z\right)^{1/2}, \\\nonumber\text{if} \quad (x_0, y_0) \quad\text{satisfies second condition of \textbf{property$^*$}},
\end{eqnarray}
where $C>0$ is a constant depending on the diameter of $\omega$ and the distance between point $(x_0, y_0)$ and the boundary of $\omega$.
\end{lemma}

\bigskip
In view of Lemma \ref{curve} one guarantees the existence of constant $C'$ depending on the distance between the boundary of $\mathbb{B}_\epsilon$ and the boundary of $\Omega$ such that 
for any $(x, y)\in \mathbb{B}_\epsilon$
\begin{eqnarray*}|u(x, y)|^2\le C' \left(\int_{l(x)}\left|\frac{\partial u}{\partial y}(x, z)\right|^2\,d z+ \int_{l(x)}|u(x, z)|^2\,d z\right)\\+C' \left( \int_{h(y)}\left|\frac{\partial u}{\partial z}(z, y)\right|^2\,d z+ \int_{h(y)}|u(z, y)|^2\,d z\right).
\end{eqnarray*}

Therefore we obtain
\begin{eqnarray*}
\int_{\mathbb{B}_\epsilon}|u(x, y)|^2\,d x\,d y= C' \int_{\mathbb{B}_\epsilon}\int_{l(x)}\left(\left|\frac{\partial u}{\partial z}(x, z)\right|^2+|u(x, z)|^2\right)\,d z\,d x\,d y\\\nonumber+ C' \int_{\mathbb{B}_\epsilon} \int_{h(y)}\left(\left|\frac{\partial u}{\partial z}(z, y)\right|^2+|u(z, y)|^2\right)\,d z\,d x\,d y\\
\le C'\int_{-\epsilon}^\epsilon\int_{-\sqrt{\epsilon^2-x^2}}^{\sqrt{\epsilon^2-x^2}}\int_{l(x)}\left(\left|\frac{\partial u}{\partial z}(x, z)\right|^2+|u(x, z)|^2\right)\,d z\,d x\,d y\\+ C' \int_{-\epsilon}^\epsilon\int_{-\sqrt{\epsilon^2-y^2}}^{\sqrt{\epsilon^2-y^2}}\int_{h(y)}\left(\left|\frac{\partial u}{\partial z}(z, y)\right|^2+|u(z, y)|^2\right)\,d z\,d y\,d x\\\le
C'\epsilon  \int_\Omega(|\nabla u|^2+ 2|u|^2)\,d x\,d y\le 2C'\epsilon \|u\|_1^2.
\end{eqnarray*} 

In view of the fact that $\epsilon\le 2\varepsilon$ the above inequality means
\begin{equation}\label{vepsilon} \int_{\mathbb B_\epsilon} |u|^2 \,d x\, d y \le  4C'\varepsilon \|u\|_1^2,
\end{equation}
which estimates the first term in the right-hand side of (\ref{ineq.}).

\bigskip
Let us now investigate the second term in the right-hand side of (\ref{ineq.}). We are going to establish the upper bound for the larger integral $\int_{\mathbb B_\epsilon}|J_1'u|^2\, d x\,d y$.

Passing to polar coordinates one gets
\begin{eqnarray*}
\nonumber &&\int_{\mathbb B_\epsilon}|J_1'u|^2\, d x\,d y = \int_0^\epsilon \int_0^{2\pi} r \left|\frac{r}{\epsilon}
\tilde{u}(\epsilon, \varphi)\right|^2\, d r\,d \varphi\\&&\le\int_0^\epsilon \int_0^{2\pi} r |\tilde{u}(\epsilon, \varphi)|^2\, d r\,d \varphi\le \epsilon \int_0^\epsilon \int_0^{2\pi} 
|\tilde{u}(\epsilon, \varphi)|^2\, d r\,d \varphi\\&&= \epsilon^2 \int_0^{2\pi} 
|\tilde{u}(\epsilon, \varphi)|^2\,d \varphi.
\end{eqnarray*}
Taking into account that $\epsilon \int_0^{2\pi} 
|\tilde{u}(\epsilon, \varphi)|^2\, d \varphi$ coincides with the curvilinear integral 
$\int_{\partial \mathbb{B}_\epsilon} |u|_{\partial \mathbb{B}_\epsilon}^2\, d \mu$, where $|u|_{\partial \mathbb{B}_\epsilon}$ is the trace of $|u|$ on the circle $\partial \mathbb{B}_\epsilon$, the above bound performs to
\begin{equation}\label{u_gamma}
\int_{\mathbb B_\epsilon}|J_1'u|^2\, d x\,d y \le \epsilon
\int_{\partial \mathbb{B}_\epsilon} |u|_{\partial \mathbb{B}_\epsilon}^2\, d \mu.
\end{equation}

The suitable bound for the above estimate is guaranteed by the following trace inequality:

\begin{lemma}\label{trace}
Let $v\in \mathcal{H}^1(\Omega\setminus \mathbb{B})$, where $\mathbb{B}$ is a ball with center in zero and radius $\tau>0$.
Then there exists a constant $C''>0$ depending on the diameter of $\Omega$ and the distance between the boundary of $\mathbb{B}$ and the boundary of $\Omega$ such that
$$\int_{\partial \mathbb B} 
|v|_{\partial \mathbb B}^2\, d \mu\le C'' \int_\Omega(|\nabla v|^2+ |v|^2)\,d x\,d y.$$
\end{lemma}

Hence, having the fact that $\epsilon\le 2\varepsilon$ inequality (\ref{u_gamma}) means
\begin{equation}\label{B}\int_{\mathbb B_\epsilon}|J_1'u|^2\, d x\,d y\le 2C'' \varepsilon \|u\|_1^2.\end{equation}

In view of the above bound and  (\ref{start 6}), (\ref{ineq.}), (\ref{vepsilon}) it follows from that
$$
\|J' u- J_1'u\|_0^2\le  8C'\varepsilon \|u\|_1^2+ 4 C'' \varepsilon \|u\|_1^2. 
$$
which proves (\ref{6}) with $\delta=\mathcal{O}(\sqrt{\varepsilon})$.\\
\\

We now give the proof of the estimate (\ref{7}), i.e. {\it under the assumptions stated in Theorem \ref{Main} inequality (\ref{7}) takes place with $k=2$ and $\delta=\mathcal{O}(\varepsilon^{1/6})$ for small enough 
$\varepsilon$.}

In view of our construction
\begin{eqnarray}\nonumber&&
|a(f, J_1'u)- a'(J_1f, u)|=\left|\int_\Omega \nabla f\,\overline{\nabla J_1' u}\,d x\,d y- \int_\Omega\nabla f\overline{\nabla u}\,d x\,d y\right|\\\label{7'}
&&= 
\left|\int_\Omega \nabla f\,\overline{(\nabla J_1' u- \nabla u)}\,d x\,d y\right|.
\end{eqnarray}
 
We are going to give the upper estimate of the right-hand side of (\ref{7'}). As before without loss of generality suppose that $u\in C^\infty(\Omega_{K_\varepsilon})$.

Using the fact that $J_1' u= u$ on $\Omega_\epsilon$ writing
\begin{equation}\label{proof.7}
\left|\int_\Omega \nabla f\,\overline{(\nabla J_1' u- \nabla u)}\,d x\,d y\right|=\left|\int_{\mathbb B_\epsilon} \nabla f\,\overline{(\nabla J_1'u- \nabla u)}\,d x\, d y\right|.
\end{equation}

In view of the Cauchy inequality, the right-hand side of (\ref{proof.7}) can be estimated as follows

\begin{eqnarray}\nonumber&&\left|\int_{\mathbb B_\epsilon} \nabla f \,\overline{(\nabla J_1'u- \nabla u)}\,d x\, d y\right|\\\nonumber&&\le\left(\int_{\mathbb B_\epsilon} |\nabla f|^2 \,d x\, d y\right)^{1/2} \,\left(\int_{\mathbb B_\epsilon}|\nabla J_1' u- \nabla u|^2\,d x\, d y\right)^{1/2}
\\\label{first. inequality}
&&\le\sqrt{2}  \left(\int_{\mathbb B_\epsilon} |\nabla f|^2 \,d x\, d y\right)^{1/2} \,\left(\int_{\mathbb B_\epsilon}|\nabla J_1' u|^2\,d x\,d y+ \int_{\mathbb B_\epsilon}|\nabla u|^2\,d x\, d y\right)^{1/2}.
\end{eqnarray}
The suitable upper bound for $\int_{\mathbb B_\epsilon}|\nabla f|^2\,d x\,d y$ is guaranteed by
the following lemma (the proof is given in the Appendix):
\begin{lemma}\label{Galpha0}
For any function $g\in \mathrm{Dom}(-\Delta^\Omega_N)$ the following estimate takes place
\[
\int_{\mathbb B_\epsilon}|\nabla g|^2 \,d x\, d y\le \tilde{C} \epsilon^{4/3}\, \|g\|^2_2,
\]
with the constant $\tilde{C}$ depends on the distance of $\mathbb{B}_\varepsilon$ from the boundary of $\Omega$.
\end{lemma}
To proceed further we need to estimate the integral
$\int_{\mathbb B_\epsilon} |\nabla J_1' u|^2\,d x\,d y$.
Passing to polar coordinates one gets
\begin{eqnarray}\nonumber&&\int_{\mathbb B_\epsilon} |\nabla J_1' u|^2\,d x\,d y=
\int_{\mathbb B_\epsilon} \left(\left|\frac{\partial (J_1' u)}{\partial x}\right|^2+ \left|\frac{\partial (J_1' u)}{\partial y}\right|^2 \right)\,d x\,d y\\\nonumber&&
= \int_0^\epsilon \int_0^{2\pi}
r \left|\frac{1}{\epsilon} \tilde{u}(\epsilon, \varphi)\,\cos \varphi -\frac{1}{\epsilon} \frac{\partial \tilde{u}}{\partial \varphi}(\epsilon, \varphi)\, \sin \varphi\right|^2\,d r\,d \varphi\\\nonumber&&+ \int_0^\varepsilon \int_0^{2\pi}
r \left|\frac{1}{\epsilon} \tilde{u}(\epsilon, \varphi)\,\sin \varphi +\frac{1}{\epsilon} \frac{\partial \tilde{u}}{\partial \varphi}(\epsilon, \varphi)\, \cos \varphi\right|^2\,d r\,d \varphi\\\nonumber&& \le \frac{4}{\epsilon} \int_0^\epsilon \int_0^{2\pi} \left(|\tilde{u}(\epsilon, \varphi)|^2+\left|\frac{\partial \tilde{u}}{\partial \varphi}(\epsilon, \varphi)\right|^2\right)\,d r\,d\varphi\\\label{J 1'}&& =4 \int_0^{2\pi}|\tilde{u}(\epsilon, \varphi)|^2\,d \varphi+ 4 \int_0^{2\pi}\left|\frac{\partial \tilde{u}}{\partial \varphi}(\epsilon, \varphi)\right|^2\,d \varphi.
\end{eqnarray}

As in proof of (\ref{6}) we again note that
$$4 \int_0^{2\pi}|\tilde{u}(\epsilon, \varphi)|^2\,d \varphi=\frac{4}{\epsilon}\,\int_{\partial \mathbb B_\epsilon} 
|u|_{\partial \mathbb B_\epsilon}^2\, d \mu,$$
which can be estimated using Lemma \ref{trace}:

\begin{equation}\label{first term}
4 \int_0^{2\pi}|\tilde{u}(\epsilon, \varphi)|^2\,d \varphi\le \frac{4 C''}{\epsilon} \|u\|_1^2.
\end{equation}

\bigskip

The upper bound for the second term of (\ref{J 1'}) can be found by the following auxiliary lemma (the proof is given in the Appendix).

\begin{lemma}\label{auxiliary}
Let $\mathbb B_{2\varepsilon}$ and $\mathbb B_\varepsilon$ be the balls with center in origin and radii 
$\varepsilon$ and $2\varepsilon$. Let $g\in \mathcal{H}^1(\mathbb B_{2\varepsilon}\setminus \mathbb B_\varepsilon)$. Then there exists $\tau\in (\varepsilon, 2\varepsilon)$ such that
\begin{equation}\label{auxiliary1}\int_0^{2\pi} \left|\frac{\partial \tilde{g}}{\partial \varphi}(\tau, \varphi)\right|^2\,d \varphi \le 4 \|g\|_1^2,\end{equation}
where $\tilde{g}(r, \varphi):= g(r\cos\varphi, r\sin\varphi)$.
\end{lemma}

Now we are able to choose the number $\epsilon$. Let us apply Lemma \ref{auxiliary} for $g=u$. In case if mentioned in above lemma $\tau$ belongs to interval 
$(\varepsilon, 3\varepsilon/2]$ then we choose $\epsilon$ as the supremum of all such numbers in
$(\varepsilon, 3\varepsilon/2]$. In the opposite case if $\tau\in (3\varepsilon/2, 2\varepsilon)$ then let $\epsilon$ be the infimum of such numbers. Since $u$ is smooth function then inequality (\ref{auxiliary1}) is satisfied with $\tau=\epsilon$. 
Combining (\ref{first term}) together with Lemma \ref{auxiliary} inequality (\ref{J 1'}) implies
\begin{equation}\label{J1'}
\int_{\mathbb B_{\epsilon}} |\nabla J_1' u|^2\,d x\,d y\le\frac{4 C''}{\epsilon}\, \|u\|_1^2+ 16 \|u\|_1^2.
\end{equation}

Finally, employing (\ref{J1'}) and Lemma \ref{Galpha0} in inequality (\ref{first. inequality}) and having that $\int_{\mathbb B_{\epsilon}}  |\nabla u|^2\,d x\, d y\le \|u\|_1^2$ it follows that
\begin{eqnarray*}&&
\left|\int_{\mathbb B_{\epsilon}} \nabla f \,\overline{(\nabla J_1'u- \nabla u)}\,d x\, d y\right|\le (2 \tilde{C})^{1/2}{\epsilon}^{2/3} \left(\frac{4 C''}{\epsilon}+ 17\right)^{1/2} \|f\|_2\,\|u\|_1 \\&&=(2 \tilde{C} (4 C''+ 17{\epsilon}))^{1/2}{\epsilon}^{1/6}\,\|f\|_2\,\|u\|_1.
\end{eqnarray*}

Using the above estimate and the fact that ${\epsilon}\le 2\varepsilon$ inequality (\ref{proof.7}) 
implies
\begin{equation}\label{as..}
\left|\int_{\Omega} \nabla f\,\overline{(\nabla J_1' u- \nabla u)}\,d x\,d y\right|\le2^{2/3}(\tilde{C} (4 C''+ 34\varepsilon))^{1/2} \varepsilon^{1/6}  \|f\|_2\,\|u\|_1.
\end{equation}

\bigskip
Returning to (\ref{7'}) and applying the bound (\ref{as..}) we have
\begin{equation}\label{as.}|a(f, J_1'u)- a'(J_1f, u)|\le 2^{2/3}(\tilde{C} (4 C''+ 34\varepsilon))^{1/2}  \varepsilon^{1/6}  \|f\|_2\,\|u\|_1.
\end{equation}

It is easy to notice that the right-hand side of inequality (\ref{as.}) for small enough $\varepsilon$ satisfies
$$
\mathrm{r.h.s.}(\ref{as.})= \mathcal{O}\left(\varepsilon^{1/6}\right)\,\|f\|_2\, \|u\|_1,
$$
which ends the proof.\qed
%%%%%%%%

\section{Appendix}

In this section we give the proofs of Lemmata \ref{curve},\, \ref{trace},\, \ref{Galpha0}  and  \ref{auxiliary}.

\subsection{Proof of Lemma \ref{curve}}
\setcounter{equation}{0}
Let $(x_0, y_0)\in \Omega_{K_\varepsilon}$. Assume the validity of the first condition of \textbf{property$^*$}. Let $y_1(x_0)$ be a point of intersection of $l(x_0)$ and the boundary of $\Omega$. Without less of generation suppose that $y_1(x_0)< y_0$. 
One can easily check that for almost any
 $(x_0, y_0)\in \omega$ there exists $y_2(x_0)\in (y_1(x_0), y_0)$ such that
$$|u(x_0, y_0)|\le\frac{1}{\sqrt{y_0-y_1(x_0)}}\sqrt{\int_{y_1(x_0)}^{y_0}|u(x_0, z)|^2\,d z}.$$

Therefore 
\begin{eqnarray*}
|u(x_0, y_0)|^2= \left|u(x_0, y_2(x_0))+\int_{y_2(x_0)}^{y_0}\frac{\partial u}{\partial z}(x_0, z)\,d z\right|^2\\\nonumber\le 2|u(x_0, y_2(x_0))|^2+ 2(y_0-y_2(x_0))\int_{y_2(x)}^y\left|\frac{\partial u}{\partial z}(x, z)\right|^2\,d y\\\nonumber\le \frac{2}{y_0-y_1(x_0)}\int_{y_1(x_0)}^{y_0}|u(x_0, z)|^2\,d z
+ 2(y_0-y_2(x_0))\int_{y_2(x_0)}^{y_0}\left|\frac{\partial u}{\partial z}(x_0, z)\right|^2\,d z\\\le 
\frac{2}{\mathrm{dist}((x_0, y_0), \partial \Omega)}\int_{y_1(x_0)}^{y_0}|u(x_0, z)|^2\,d z
+ 2 \mathrm{diam}(\Omega)\int_{y_1(x_0)}^{y_0}\left|\frac{\partial u}{\partial z}(x_0, z)\right|^2\,d z\\\le\frac{2}{\mathrm{dist}((x_0, y_0), \partial \Omega)}\int_{l(x_0)}|u(x_0, z)|^2\,d z
+ 2 \mathrm{diam}(\Omega)\int_{l(x_0)}\left|\frac{\partial u}{\partial z}(x_0, z)\right|^2\,d z,
\end{eqnarray*} 
where  $\mathrm{diam}(\Omega)$ is the diameter of $\Omega$ and $\mathrm{dist}((x_0, y_0), \partial\Omega)$ is the distance between $(x_0, y_0)$ and the boundary of $\Omega$. This proves (\ref{curve1}).

The case when $(x_0, y_0)$ satisfies second condition of \textbf{property$^*$} can be studied similarly. Repeating the same ideas one gets the validity of (\ref{curve2}). \qed

\subsection{Proof of Lemma\,\ref{trace}}
\setcounter{equation}{0}

Let us estimate separately the curvilinear integrals $\int_{\gamma^i}|v|_{\partial \gamma^i}^2\,d \mu,\,i=1,2,3,4$, where $\gamma^i$ are given in polar coordinates as follows 
\begin{eqnarray*}
\gamma^1=\left\{r=\tau, \frac{\pi}{4}\le\varphi\le\frac{3\pi}{4}\right\},  \,\,
\gamma^2=\left\{r=\tau, \frac{3\pi}{4}\le\varphi\le\frac{5\pi}{4}\right\},\\
\gamma^3=\left\{r=\tau, \frac{5\pi}{4}\le\varphi\le\frac{7\pi}{4}\right\},\,\,
\gamma^4=\left\{r=\tau, -\frac{\pi}{4}\le\varphi\le\frac{\pi}{4}\right\}
\end{eqnarray*}
and $|v|_{\gamma^i}$ is the trace of $v$ on $\gamma^i$.

Let us start from $\int_{\gamma^1}|v|_{\gamma^1}^2\,d \mu$. We have 
\begin{eqnarray*}
\int_{\gamma^1} |v|_{\gamma^1}^2\,d \mu=\tau \int_{-\tau/\sqrt{2}}^{\tau/\sqrt{2}}\frac{1}{\sqrt{\tau^2-x^2}}|v(x, \sqrt{\tau^2-x^2})|^2\,d x\\\nonumber\le \sqrt{2} \int_{-\tau/\sqrt{2}}^{\tau/\sqrt{2}}|v(x, \sqrt{\tau^2-x^2})|^2\,d x.
\end{eqnarray*}

Since $(x, \sqrt{\tau^2-x^2}))$ for each $x, |x|\le \tau/\sqrt{2}$, obviously satisfies first condition of \textbf{property$^*$} for domain $\Omega_{K_\varepsilon}$ then using estimate \ref{curve1} one gets
\begin{eqnarray}
\nonumber\int_{\gamma^1} |v|_{\gamma^1}^2\,d \mu \le \sqrt{2}C^2 \int_{-\tau/\sqrt{2}}^{\tau/\sqrt{2}} \left(\int_{l(x)}\left|\frac{\partial v}{\partial z}(x, z)\right|^2\,d z+ \int_{l(x)}|v(x, z)|^2\,d z \right)\,d x\,d y\\\nonumber\le
\sqrt{2} C^2  \int_{\Omega\cap\{|x|\le\tau/\sqrt{2}\}}|v(x, y)|^2\,d x\,d y\\\nonumber
+\sqrt{2} C^2 \int_{\Omega\cap\{|x|\le\tau/\sqrt{2}\}}\left|\frac{\partial v}{\partial y}(x, y)\right|^2\,d x\,d y\\\label{lastgamma1}\le
\sqrt{2} C^2 \int_\Omega|v(x, y)|^2\,d x\,d y
+ \sqrt{2} C^2\int_\Omega\left|\frac{\partial v}{\partial y}(x, y)\right|^2\,d x\,d y.
\end{eqnarray}

In the same way we are able to estimate the curvilinear integral 
\begin{eqnarray}\label{gamma3}\int_{\gamma^3}|v|_{\gamma^3}^2\,d \mu\le \sqrt{2} C^2 \int_\Omega|v(x, y)|^2\,d x\,d y\\\nonumber
+ \sqrt{2} C^2 \int_\Omega\left|\frac{\partial v}{\partial y}(x, y)\right|^2\,d x\,d y.
\end{eqnarray}

In a similar way one deals with the integrals over curves $\gamma^2, \gamma^4$. We establish
\begin{eqnarray*}
\int_{\gamma^2}|v|_{\gamma^2}^2\,d \mu\le \sqrt{2} C^2 \int_\Omega|v(x, y)|^2\,d x\,d y\\
+ \sqrt{2} C^2 \int_\Omega\left|\frac{\partial v}{\partial x}(x, y)\right|^2\,d x\,d y,\\
\int_{\gamma^4}|v|_{\gamma^4}^2\,d \mu\le \sqrt{2} C^2\int_\Omega|v(x, y)|^2\,d x\,d y\\
+ \sqrt{2} C^2 \int_\Omega\left|\frac{\partial v}{\partial x}(x, y)\right|^2\,d x\,d y.
\end{eqnarray*}

Combining the above bounds together with (\ref{lastgamma1}) and (\ref{gamma3}) we finish the proof of lemma.
\qed

\subsection{Proof of Lemma \ref{Galpha0}}
\setcounter{equation}{0}

To proceed with a proof we need the following auxiliary material  \cite[Lemma\,4.9]{MK06}:
\begin{lemma}
Let $\Pi\subset \mathbb{R}^n$ be a convex set and let $G$ and $Q$ be arbitrary measurable sets in $\Pi$ with $\mu\, (G)\ne 0$.
Then, for all $v\in \mathcal{H}^1(\Pi)$, the following inequality holds:
\begin{eqnarray}\label{Marchenko'}
&&\int_Q |v|^2\,d x\,d y\\ \nonumber
&&\le \frac{2\mu\, (Q)}{\mu \,(G)} \int_G |v|^2\,d x\,d y
+\frac{C(n) (d(\Pi))^{n+1} (\mu\, (Q))^{1/n}}{\mu \,(G)} \int_\Pi |\nabla v|^2\, d x\,d y,
\end{eqnarray}  
where $d(\Pi)$ is the parameter of $\Pi$, $\mu$ is the Lebesque measure on $\mathbb{R}^n$, and the constant $C(n)$ depends only on the dimension of $\mathbb{R}^n$.
\end{lemma}

Let $G_\epsilon$ be a convex subset of $\Omega$ to be chosen later which contains $\mathbb B_\epsilon$. 
Since $g\in \mathcal{H}^2_{\mathrm{loc}}(\Omega)$ in view of the interior regularity theory then applying (\ref{Marchenko'}) for $Q=\mathbb B_\epsilon$ and $G=\Pi= G_\epsilon$ one infers that
\begin{eqnarray*}&& \int_{\mathbb B_\epsilon} \left|\frac{\partial g}{\partial x}\right|^2 \,d x\, d y\le 
\frac{2 \mu \,(\mathbb B_\epsilon)}{\mu\, (G_\epsilon)} \int_{G_\epsilon} \left|\frac{\partial g}{\partial x}\right|^2\,d x\,d y\\&&
+\frac{C(2) (d\,(G_\epsilon))^3 (\mu\, (\mathbb B_\epsilon))^{1/2}}{\mu \,(G_\epsilon)} \int_{G_\epsilon}\left|\nabla\left(\frac{\partial g}{\partial x}\right)\right|^2\,d x\,d y
\\
&&=\frac{2 \mu \,(\mathbb B_\epsilon)}{\mu\, (G_\epsilon)} \int_{G_\epsilon} \left|\frac{\partial g}{\partial x}\right|^2\,d x\,d y+
 \frac{C(2) (d\,(G_\epsilon))^3 (\mu\, (\mathbb B_\epsilon))^{1/2}}{\mu \,(G_\epsilon)} \int_{G_\epsilon}\left(\left|\frac{\partial^2 g}{\partial x^2}\right|^2+ \left|\frac{\partial^2 g}{\partial x \partial y}\right|^2\right)\,d x\,d y\end{eqnarray*}
 and
 \begin{eqnarray*}&&
 \int_{\mathbb B_\epsilon} \left|\frac{\partial g}{\partial y}\right|^2 \,d x\, d y \le \frac{2 \mu \,(\mathbb B_\epsilon)}{\mu\, (G_\epsilon)} \int_{G_\epsilon} \left|\frac{\partial g}{\partial y}\right|^2\,d x\,d y\\&&+
 \frac{C(2) (d\,(G_\epsilon))^3 (\mu\, (\mathbb B_\epsilon))^{1/2}}{\mu \,(G_\epsilon)} \int_{G_\epsilon}\left|\nabla\left(\frac{\partial g}{\partial y}\right)\right|^2\,d x\,d y
\\&&=\frac{2 \mu \,(\mathbb B_\epsilon)}{\mu\, (G_\epsilon)} \int_{G_\epsilon} \left|\frac{\partial g}{\partial y}\right|^2\,d x\,d y+
 \frac{C(2) (d\,(G_\epsilon))^3 (\mu\, (\mathbb B_\epsilon))^{1/2}}{\mu \,(G_\epsilon)} \int_{G_\epsilon}\left(\left|\frac{\partial^2 g}{\partial y \partial x}\right|^2+ \left|\frac{\partial^2 g}{\partial y^2}\right|^2\right)\,d x\,d y.
\end{eqnarray*}

Combining the above inequalities we arrive

\begin{eqnarray}\nonumber&&
\int_{\mathbb B_\epsilon}|\nabla g|^2 \,d x\, d y\le \frac{2 \mu \,(\mathbb B_\epsilon)}{\mu\, (G_\epsilon)} \int_{G_\epsilon} |\nabla g|^2\,d x\,d y \\\nonumber&&+
 \frac{C(2) (d\,(G_\epsilon))^3 (\mu\, (\mathbb B_\epsilon))^{1/2}}{\mu \,(G_\epsilon)} \int_{G_\epsilon}\left(\left|\frac{\partial^2 g}{\partial x^2}\right|^2+ 2\left|\frac{\partial^2 g}{\partial x \partial y}\right|^2+ \left|\frac{\partial^2 f}{\partial y^2}\right|^2\right)\,d x\,d y\\\nonumber&& \le  \frac{(\mu\,(\mathbb B_\epsilon))^{1/2}}{\mu\, (G_\epsilon)}\left(2(\mu\,(\mathbb B_\epsilon))^{1/2}+ C(2) (d\,(G_\epsilon))^3 \right)\|g\|_{\mathcal{H}^2(G_\epsilon)}^2
 \\\label{nabla.f}
&&= \frac{\sqrt{\pi} \epsilon }{\mu\, (G_\epsilon)}\left(2 \sqrt{\pi} \epsilon+ C(2) (d\,(G_\epsilon))^3 \right)\|g\|_{\mathcal{H}^2(G_\epsilon)}^2.
\end{eqnarray} 

Let us choose $G_\epsilon=\{r:\, 0\le r\le\epsilon^{1-\alpha}\}$ with some $\alpha\in \left(\frac{1}{2}, 1\right)$ to be chosen later. One can easily notice that $\mathbb B_\epsilon \subset G_\epsilon$. The inequality (\ref{nabla.f}) performs to

\begin{equation}\label{Gepsilon}
\int_{\mathbb B_\epsilon}|\nabla g|^2 \,d x\, d y\le 2\epsilon 
\left(\epsilon^{2\alpha-1}+ \frac{4 C(2)}{\sqrt{\pi}} \epsilon^{1-\alpha}\right)\, \|g\|^2_{\mathcal{H}^2(G_\epsilon)}.
\end{equation}

%%%%%%%%%%%%%%%%%%
To make the further estimates we need the interior regularity theorem \cite{B16} and the following lemma:

\begin{theorem}(Interior Regularity Theorem.)
Suppose that $h\in \mathcal{H}^1(\Omega)$ is a weak solution to $-\Delta h=w$. Then $h\in \mathcal{H}^2_{\mathrm{loc}}(\Omega)$ and for each $\Omega_0\subset  \Omega$ there is a constant $c= c(\Omega_0)$ independent of $h$ and $w$ such that: 
\begin{equation}\label{Elliptic}\|h\|_{\mathcal{H}^2(\Omega_0)}\le c \left(\|h\|_{L^2(\Omega)}+\|w\|_{L^2(\Omega)}\right).
\end{equation}
\end{theorem}
\begin{lemma}\label{Delta}[The proof will be given in the next subsection]
For any $z\in \mathrm{Dom}(-\Delta_N^\Omega)$ the following estimate is valid
$$\int_\Omega |-\Delta z + z|^2\,d x\,d y \ge \frac{1}{16} \int_\Omega (|\Delta z|^2+ |z|^2)\,d x\,d y.$$
\end{lemma}

Let us go back to inequality (\ref{Gepsilon}). Using (\ref{Elliptic}) and Lemma\,\ref{Delta} one establishes that
\begin{eqnarray}\label{FF}\nonumber&&
\int_{\mathbb B_\epsilon}|\nabla g|^2 \,d x\, d y\le 4\, c^2 \epsilon 
\left(\epsilon^{2\alpha-1}+ \frac{4 C(2)}{\sqrt{\pi}} \epsilon^{1-\alpha}\right)\, (\|\Delta g\|_{L^2(\Omega)}^2+ \|g\|_{L^2(\Omega)}^2)
\\\nonumber&&\le64\, c^2 \epsilon 
\left(\epsilon^{2\alpha-1}+ \frac{4 C(2)}{\sqrt{\pi}} \epsilon^{1-\alpha}\right)\, \|-\Delta g+ g\|_{L^2(\Omega)}^2\\\label{C''}&&=64\, c^2 \epsilon 
\left(\epsilon^{2\alpha-1}+ \frac{4 C(2)}{\sqrt{\pi}} \epsilon^{1-\alpha}\right)\, \|g\|_2^2\le \tilde{C} 
\left(\epsilon^{2\alpha}+ \epsilon^{2-\alpha}\right)\, \|g\|_2^2,
\end{eqnarray} 
where $\tilde{C}:= 64 c^2 \mathrm{max}\left\{1, \frac{4 C(2)}{\sqrt{\pi}}\right\}$ and $c$ depends on the distance of $\mathbb{B}_\varepsilon$ from the boundary of $\Omega$. 

Let us investigate the function $F(\alpha):= \epsilon^{2\alpha}+ \epsilon^{2-\alpha}$ in (\ref{FF}) on interval $\left(\frac{1}{2}, 1\right)$. It attains its minimum at $\alpha_0= \frac{2}{3}- \frac{1}{3}\frac{\ln 2}{\ln \epsilon}$ and takes the value $F(\alpha_0)=\epsilon^{4/3}\left(\frac{1}{4^{1/3}}+ 2^{1/3}\right)$. 

Let us now choose $\alpha= \alpha_0$ in inequality (\ref{C''}). Then 
$$\int_{\mathbb B_\epsilon}|\nabla g|^2 \,d x\, d y\le \epsilon^{4/3} \left(\frac{1}{4^{1/3}}+ 2^{1/3}\right)\tilde{C}\,\|g\|_2^2,$$
which concludes the proof of Lemma \ref{Galpha0} with $C'=\left(\frac{1}{4^{1/3}}+ 2^{1/3}\right)\tilde{C}$.\qed

\subsection{Proof of Lemma \ref{auxiliary}}

We first prove the following auxiliary statement: there exists $\tau\in (\varepsilon, 2\varepsilon)$ such that

\begin{equation}\label{statement}\int_0^{2\pi} \left|\nabla g(\tau \cos \varphi, \tau \sin \varphi)\right|^2\,d \varphi \le \frac{1}{\varepsilon^2} \|g\|_1^2.\end{equation}

Let us assume the opposite: for any $r\in (\varepsilon, 2\varepsilon)$ we have 
$$
\int_0^{2\pi} \left|\nabla g(r \cos \varphi, r \sin \varphi)\right|^2\,d \varphi>\frac{1}{\varepsilon^2} \|g\|_1^2.
$$

Passing to polar coordinates in integral $\int_{\mathbb B_{2\varepsilon}\backslash \mathbb B_\varepsilon} |\nabla g|^2\,d x\,d y$ and using the above bound we get
$$\int_{\mathbb B_{2\varepsilon}\backslash \mathbb B_\varepsilon} |\nabla g|^2\,d x\,d y= \int_\varepsilon^{2\varepsilon} \,\int_0^{2\pi} r  \left|\nabla g(r \cos \varphi, r \sin \varphi)\right|^2\,d \varphi\, d r>\|g\|_1^2.$$

This contradicts with the fact that the left-hand side of the above inequality does not exceed $\|g\|_1^2$. Hence there exists at least one number $\tau\in (\varepsilon, 2\varepsilon)$ such that 

\begin{equation}\label{varphi0}
\int_0^{2\pi} |\nabla g(\tau \cos \varphi, \tau \sin \varphi)|^2\,d \varphi\le \frac{1}{\varepsilon^2} \|g\|_1^2.
\end{equation}

In view of the representation of the derivative of function $\tilde{g}(r, \varphi)=g(r\cos\varphi, r\sin\varphi)$
$$\frac{\partial \widetilde{g}}{\partial\varphi}(\tau, \varphi)= -\tau \frac{\partial g}{\partial x}(\tau \cos \varphi, \tau \sin \varphi) \sin \varphi +\tau \frac{\partial g}{\partial y}(\tau \cos \varphi, \tau \sin \varphi) \cos \varphi$$ 
and the Cauchy inequality we conclude that 
$$\left|\frac{\partial \widetilde{g}}{\partial\varphi}(\tau, \varphi)\right|\le 2\varepsilon \left| \nabla g(\tau \cos \varphi, \tau \sin \varphi)\right|.
$$

Hence, employing (\ref{varphi0}) we establish that
$$
\int_0^{2\pi}\left|\frac{\partial \tilde{g}}{\partial\varphi}(\tau, \varphi)\right|^2\,d \varphi \le 4\varepsilon^2 \int_0^{2\pi}\left|\nabla g(\tau \cos \varphi, \tau \sin \varphi)\right|^2\,d \varphi\le 4\|g\|_1^2,
$$
which concludes the proof.\qed

\subsection{Proof of Lemma\,\ref{Delta}}

It is straightforward to check that
\begin{eqnarray}\nonumber\int_\Omega |-\Delta u+ u|^2\,d x\,d y=
\int_\Omega |\Delta u|^2\,d x\,d y+2 \int_\Omega |\nabla u|^2\,d x\,d y+\int_\Omega |u|^2\,d x\,d y\\\label{resolvent}\ge\int_\Omega |u|^2\,d x\,d y.\end{eqnarray}

Let us consider two cases:
\begin{eqnarray}
\label{delta1}&&\int_\Omega |\Delta u|^2\,d x\,d y \ge 4\int_\Omega |u|^2\,d x\,d y,
\\\label{delta2}&&\int_\Omega |\Delta u|^2\,d x\,d y< 4\int_\Omega |u|^2\,d x\,d y.\end{eqnarray}
Starting from the first one and employing (\ref{resolvent}) we have
\begin{eqnarray*}&&\sqrt{\int_\Omega |-\Delta u+ u|^2\,d x\,d y}\ge \sqrt{\int_\Omega |\Delta u|^2\,d x\,d y}-
\sqrt{\int_\Omega |u|^2\,d x\,d y} \\&&\ge\frac{1}{2} \sqrt{\int_\Omega |\Delta u|^2\,d x\,d y}\ge  \frac{1}{4}\sqrt{\int_\Omega |\Delta u|^2\,d x\,d y}+\frac{1}{2}\sqrt{\int_\Omega |u|^2\,d x\,d y}\\&&\ge\frac{1}{4}\left(\sqrt{\int_\Omega |\Delta u|^2\,d x\,d y}+\sqrt{\int_\Omega |u|^2\,d x\,d y}\right).
\end{eqnarray*}
Hence we arrive at the bound
\begin{equation}\label{last}\int_\Omega |-\Delta u+ u|^2\,d x\,d y\ge \frac{1}{16}\left(\int_\Omega |\Delta u|^2\,d x\,d y+\int_\Omega |u|^2\,d x\,d y\right).
\end{equation}
Now let us consider the case (\ref{delta2}). In view of inequality (\ref{resolvent}) we conclude
\begin{eqnarray*}
\int_\Omega |-\Delta u+ u|^2\,d x\,d y\ge \int_\Omega |u|^2\,d x\,d y\ge\frac{1}{2}\int_\Omega |u|^2\,d x\,d y+ \frac{1}{8} \int_\Omega |\Delta u|^2\,d x\,d y\\\ge \frac{1}{8}\left(\int_\Omega |\Delta u|^2\,d x\,d y+ \int_\Omega |u|^2\,d x\,d y\right).
\end{eqnarray*}

Combining the above estimate together with (\ref{last}) we complete the proof of the lemma. \qed

\subsection*{Acknowledgements}
D.B. acknowledges support from the Czech Science Foundation (GACR) within the project 21-07129S.

\bigskip
${}$
\\
{\bf Conflict of Interests}
\\
The authors declare that they have no conflict of interest regarding the publication of this paper.
\\
\\
{\bf Authors' contributions}
\\
All authors contributed equally to the manuscript and typed, read, and approved the final form of the manuscript, which is the result of an intensive collaboration. 
\\

\end{document}